\documentclass{jmlr}

\title[QP Based Constrained Optimization for Reliable PINN Training]{QP Based Constrained Optimization for Reliable PINN Training}

\usepackage{amsfonts}

\usepackage{mathtools}
\usepackage{siunitx}

\newtheorem{assumption}{Assumption}

\newcommand{\assums}{\ref{assum:unique_solution} - \ref{assum:gradient_condition} }

\newcommand{\Vb}{V_\beta}
\newcommand{\lr}{\gamma}
\newcommand{\param}{\theta}

\newcommand{\paramt}{\param^{(t)}} 
\newcommand{\paramtp}{\param^{(t+1)}} 

\DeclareMathOperator*{\argmin}{arg\,min}

\usepackage{times}

\author{%
 \Name{Alan Williams} \Email{awilliams@lanl.gov}\\
 \addr Applied Electrodynamics Group (AOT-AE), Los Alamos National Laboratory
 \AND
 \Name{Christopher Leon} \Email{cleon@lanl.gov}\\
 \addr Applied Electrodynamics Group (AOT-AE), Los Alamos National Laboratory
 \AND
 \Name{Alexander Scheinker} \Email{ascheink@lanl.gov}\\
 \addr Applied Electrodynamics Group (AOT-AE), Los Alamos National Laboratory
}

\begin{document}

\maketitle

\begin{abstract}
Physics-Informed Neural Networks (PINNs) have emerged as a powerful tool for integrating physics-based constraints and data to address forward and inverse problems in machine learning. Despite their potential, the implementation of PINNs are hampered by several challenges, including issues related to convergence, stability, and the design of neural networks and loss functions. In this paper, we introduce a novel training scheme that addresses these challenges by framing the training process as a constrained optimization problem. Utilizing a quadratic program (QP)-based gradient descent law, our approach simplifies the design of loss functions and guarantees convergences to optimal neural network parameters. This methodology enables dynamic balancing, over the course of training, between data-based loss and a partial differential equation (PDE) residual loss, ensuring an acceptable level of accuracy while prioritizing the minimization of PDE-based loss. We demonstrate the formulation of the constrained PINNs approach with noisy data, in the context of solving Laplace's equation in a capacitor with complex geometry. This work not only advances the capabilities of PINNs but also provides a framework for their training.
\end{abstract}


\section{Introduction}
Physics-Informed Neural Networks (PINNs) (\cite{raissi2019physics}) represent a groundbreaking approach that combines physics constraints, modeled as partial differential equation (PDE) residuals, with gathered data to solve both forward and inverse problems in machine learning (ML). With complete knowledge of the physics and boundary conditions, PINNs can be used to simply solve a PDE, with a loss function given as the sum of the PDE residual and data based losses. But there are also many situations where the PINNs optimization problem may be difficult to formulate when different parts of the loss functions ``fight'' each other. For example, in experimental contexts, fitting a neural network (NN) to potentially noisy measurements may not exactly agree with the PDE. Or, in a setting where the PDE residual may be an incomplete description of the physics, one may not want to minimize the residual to precisely zero. Network architecture and expressiveness also dictate how small losses can be made, when training. Despite their potential, PINNs face significant challenges (\cite{cuomo2022scientific}). As \cite{antonion2024machine} points out, ``Despite significant strides in augmenting PINN capabilities through published works, numerous unresolved issues persist. These encompass a broad spectrum, spanning from theoretical considerations --- such as convergence and stability --- to implementation challenges, including boundary condition management, neural network design, general PINN architecture, and optimization aspects."

We propose a training scheme aimed at simplifying the design of the loss function and ensuring convergence and stability of the training algorithm. Our approach frames the training as a constrained optimization problem, leveraging a quadratic program (QP)-based gradient descent law. Essentially, the QP-based gradient descent (QPGD) scheme dynamically balances different loss functions, shifting focus between data-based losses and PDE residual losses. The balancing is performed in a nonsmooth but continuous way, switching between the gradients of the losses. As an early case study, we consider an inverse problem to determine the voltage of a plate in a capacitor, with the potential field described by Laplace's equation, given noisy measurements. We therefore are interested in achieving an acceptable level of data based loss, determined by inherent noise level in the data, while minimizing the PDE based loss under this constraint.

A recently discovered QP based scheme gradient descent scheme, which provides the basis for this work, has been shown to have strong convergence properties (\cite{williams2022practically, williams2023experimental, williamsCDC2023, allibhoy2023control}), while being able to handle equality and inequality constraints. The design of \cite{williamsCDC2023} presents a so called ``extremum seeking'' algorithm, which studies the convergence of the constrained optimization scheme when the gradients of the problem are unknown and must be estimated online --- here a single inequality constraint is considered. The work of \cite{allibhoy2023control} studies the more general form of the nonlinear programming scheme, proving local stability under several inequality and several equality constraints. 

In the machine learning literature, several approaches exist which help mitigate issues for balancing various terms in the loss function of a PINN. A learning rate annealing (\cite{wang2021understanding}) algorithm seeks to mitigate a particular kind of instability arising from situations in which the gradients corresponding to different elements in the loss function are of different scales, and GradNorm (\cite{chen2018gradnorm}) also provides a solution to this problem. ReLoBRaLo (\cite{bischof2021multi}) develop heuristics which are used to adaptively tune the weights of the total loss function for more desirable training in the case of three benchmark PDEs. Our method, although distinct from the aforementioned work, addresses similar issues by a reformulation of the overall optimization problem and presents an algorithm which solves it.

While PINNs do not guarantee hard physics constraints, but only gently push a neural network's output towards satisfying physics via additional terms in the cost function, Physics Constrained Neural Networks (PCNNs) have been developed with hard physics constraints (\cite{scheinker2023physics}). For example, PCNNs have been developed utilizing three-dimensional convolutional neural network-based neural operators with hard physics constraints that enforce the spatiotemporal Maxwell's equations built into the network architectures for electrodynamics by generating vector and scalar potential fields from which the electromagnetic fields are created (\cite{scheinker2023physics}). One of the difficulties of the PCNN approach is that each problem is handled with a unique approach, while PINNs are much more easily applied to any physical problem by simply adding PDE-based terms to the training cost function. Our case study is motivated by the explicit presence of uncertainty. While the incorporation of uncertainty into scientific ML has been studied through a variety of techniques (\cite{garcia2024machine, rautela2024time, psaros2023uncertainty, yang2019adversarial, zhang2019quantifying}), here we take the approach that we require the network try to satisfy the known physics while encouraging the network to fit the data only up to some error. 

\section{Problem Statement and Intuition}
Consider the constrained optimization problem
\begin{equation} \label{eqn:constrained_opt_problem}
    \min_\param f(\param) \text{ subject to } \param \in \mathcal{S},
\end{equation}
where 
    \begin{equation}
        \mathcal{S} \coloneq \{\param \in \mathbb{R}^n: g(\param) \leq 0\}.
    \end{equation}
with the objective $f: \mathbb{R}^n \mapsto \mathbb{R}$ and constraint $g: \mathbb{R}^n \mapsto \mathbb{R}$. 

The algorithm presented, which solves \eqref{eqn:constrained_opt_problem}, is inspired by the QP based safe controllers first formulated by \cite{ames2016control}. To understand the algorithm, consider first the continuous time dynamics of gradient descent of the objective
\begin{equation}
    \dot \param (t) = - \nabla f(\param) = u(\param) \label{eqn:simple_sys}
\end{equation}
where $u(\param)$ is a ``nominal" controller solving the unconstrained problem -- simply minimizing $f$. In this approach, one can interpret $g$ as a control barrier function (CBF). In particular, $g$ is a CBF-like function (of the ``zeroing" type (\cite{ames2016control})) for (\ref{eqn:simple_sys}) with respect to $\mathcal{S}$. (Note that in the work of \cite{ames2016control}, they use $h$ to define the constraint and consider feasible points as $h(\param)\geq0$.) The equivalent conditions below, if they hold, 
\begin{equation}
    \nabla g(\param)^T u(\param) + c g(\param) \leq 0 \iff  \dot{g} +c g(\param) \leq 0 \label{eqn:bar_condition}
\end{equation}
imply two important properties: i) $\mathcal{S}$ is forward invariant because $g(\param) = 0 \implies \dot g \leq 0$ and ii) trajectories are attracted to $\mathcal{S}$ because $g(\param) >0 \implies \dot g < 0$. Of course this condition is not automatically guaranteed by $u(\param)$ and so one can formulate the QP below to design a control law $u_s = u(\param) + \bar{u}$ satisfying
\begin{align}
    \bar{u} = \argmin_{v \in \mathbb{R}^n} ||v||^2 \quad \text{subject to} \quad \nabla g(\param)^T (u(\param) + v) + c g(\param) \leq 0,
\end{align}
which finds the smallest additional action to add to the nominal law, to guarantee \eqref{eqn:bar_condition} with respect to $\dot \theta = u_s$. If $||\nabla g(\param) || \neq 0$, an explicit solution for $\bar{u}$ exists:
\begin{equation}
\bar{u} =  \dfrac{\max\{ - \nabla f(\param)^T \nabla g(\param) + c g(\param),0\}}{||\nabla g(\param) ||^2} \nabla g(\param). \label{eqn:sf_term_exact}
\end{equation}

\section{Design and Asymptotic Convergence}
QPGD is based upon the discretization of the parameter update law $\dot \theta = u(\param) + \bar u$ for $\bar u$ in \eqref{eqn:sf_term_exact}, where we take a small ``learning rate'' as the discretized change in time. With $t \in \mathbb{N}$, learning rate $\lr>0$, the discrete update law is
\begin{equation} \label{eqn:update_law}
    \paramtp = \paramt - \lr \left( \nabla f(\paramt) + \alpha(\paramt) \nabla g(\paramt) \right),
\end{equation}
where $\alpha(\paramt) \geq 0$ is
\begin{align}
    \alpha(\paramt)  &= 
    \displaystyle\frac{\max\{-\nabla f(\paramt)^T \nabla g(\paramt) + c g(\paramt), 0\}}{\max\{||\nabla g(\paramt) ||^2, \epsilon_\alpha \}} \label{eqn:alpha} 
\end{align}
and $\nabla$ represents the gradient, $|| \cdot ||$ is the $l^2$-norm, and $\epsilon_\alpha>0$ is a small constant, avoiding a division by zero. The parameter $c > 0$ governs the rate of the approach to the constraint boundary.

Consider the following set of assumptions:
\begin{assumption} \label{assum:unique_solution}
    There exists a unique $\param^*$ which solves \eqref{eqn:constrained_opt_problem} and if $\nabla f(\param) = 0$ for $\param \in \mathcal{S}$, then $\param = \param^*$.
\end{assumption}

\begin{assumption} \label{assum:fg_assums}
    $f, g : \mathbb{R}^n \to \mathbb{R}$ are differentiable with Lipschitz continuous gradients, having Lipschitz constants $L_f, L_g > 0$.
\end{assumption}

\begin{assumption} \label{assum:safe_set}
    The function is $g$ radially unbounded and there exists an $a \in \mathbb{R}^n$ such that $g(a) < 0$.
\end{assumption}

\begin{assumption} \label{assum:gradient_condition}
    The gradient $\nabla g(\param) : \mathbb{R}^n \to \mathbb{R}^n$ does not vanish for values of $x$ in or near the region defined by $g(\param)\geq 0$. Namely, there exists a pair $\epsilon_g, l_g > 0$ such that
\begin{equation}
    || \nabla g(\param) || \geq l_g \text{ for all } \param \in \{y : g(y) \geq - \epsilon_g \}.
\end{equation}
\end{assumption}

For simplicity of the presentation, Assumption \ref{assum:unique_solution} precludes the possibility of i) studying a set of minimizers which solve \eqref{eqn:constrained_opt_problem} ii) saddle points of $f$ existing on $\mathcal{S}$. With many local minimizers the scheme in \eqref{eqn:update_law} will converge locally and not globally. It is also possible to study an unbounded set $\mathcal{S}$, and slightly relax Assumption \ref{assum:safe_set} (see \cite{williams2024semiglobal}).

The main theoretical result is asymptotic convergence to the optimum from any initial parameter, in a semiglobal sense with respect to the learning rate.
\begin{theorem} \label{}
    Let Assumptions \assums hold. There exists an $\epsilon_\alpha^*$ such that for any $\epsilon_\alpha \in (0, \epsilon_\alpha^*)$ and for any initial parameter $\param^{(0)} \in \mathbb{R}^n$ there exists a learning rate $\lr^* > 0$ such that for all $\lr \in (0, \lr^*)$ the sequence $\{ \paramt \}$ asymptotically converges to $\param^*$ as $t\to \infty$ given the the update law $\eqref{eqn:update_law}$.
\end{theorem}
See appendix for proof sketch. The key idea in proving convergence is the use of a Lyapunov function of the form
\begin{equation} \label{eqn:Vb_definition}
    \Vb(\param) = f(\param) - f^* + \beta \max\{g(\param), 0\}. 
\end{equation}
This function has a minimum value of $0$, at the optimum $\param^*$, for a sufficiently large $\beta$, under the assumptions provided. In optimization literature this function is often referred to as a ``penalty function'' with recent use in analysis of safe extremum seeking by \cite{williams2023experimental, williams2024semiglobal} and gradient flow systems by \cite{allibhoy2023control}. 

\textbf{On implementation:} The popular optimizer Adam (\cite{kingma2014adam}) is commonly used to speed up training through the use of momentum and RMSprop. When we compute the term $ \nabla f(\paramt) + \alpha(\paramt) \nabla g(\paramt)$ in \eqref{eqn:update_law}, we allow Adam, with momentum and RMSProp, to use this ``mixed'' gradient in order to update the internal exponential moving averages and perform the update to the parameter -- this, in addition to a learning rate schedule, speeds up our training in practice. 

In the next section, the fully connected neural network and the loss function, likely do not satisfy our assumptions and therefore there may be many solutions, or many unconnected islands of local solutions, and the level set $\{g(\param) > \bar g \}$ may be a collection of disconnected sets (Assumption \ref{assum:gradient_condition} precludes ``islands''). QPGD guarantees approximate invariance of the level $\{g(\param) > \bar g \}$ for $\bar g > 0$. Therefore, it may be desirable to first pretrain the PINNs with a static loss $f(\param) + \lambda g(\param)$ to move the parameters close to a particular ``island'' of $\{g(\param) > \bar g \}$, before implementing \eqref{eqn:update_law}. It is also possible to clip the value of $\alpha(\param)$ at some positive number, especially if one does not pretrain. The convergence of \eqref{eqn:update_law} may also be guaranteed with additional assumptions when $\alpha$ is clipped, but for simplicity we do not include it in our presentation in this section. If the PINN is not pretrained, and $\alpha$ is not clipped, a smaller learning rate may be required for stability.

\section{Case Study: Capacitor with Complex Geometry}
PINNs can be adapted to find solutions to real physical systems (\cite{raissi2020hidden, cai2021physics, jarolim2023probing}) and applied to solving inverse problems, something traditional PDE solvers struggle with.  A straightforward formulation would include a loss function with a term corresponding to satisfying the PDE, boundary conditions (partial or full), and a term for fitting the data: $l = l_{\text{PDE}} + l_{\text{BC}} + l_{\text{DATA}}$. One issue with this approach, however, is that given measurement uncertainty, the network is likely being asked to solve an impossible task: satisfy both the known PDE and the noise corrupted data. Given the competing terms in the loss function, the network may sacrifice satisfying the known physics to fit the inaccurate data, possibly leading to overfitting. One may try to solve this by giving greater weight to the PDE term in the loss function, but this suffers from the question of how to choose the weighing hyperparameter and offers no guarantees. 

\begin{figure}[t]
\centering  \includegraphics[width=0.99\columnwidth]{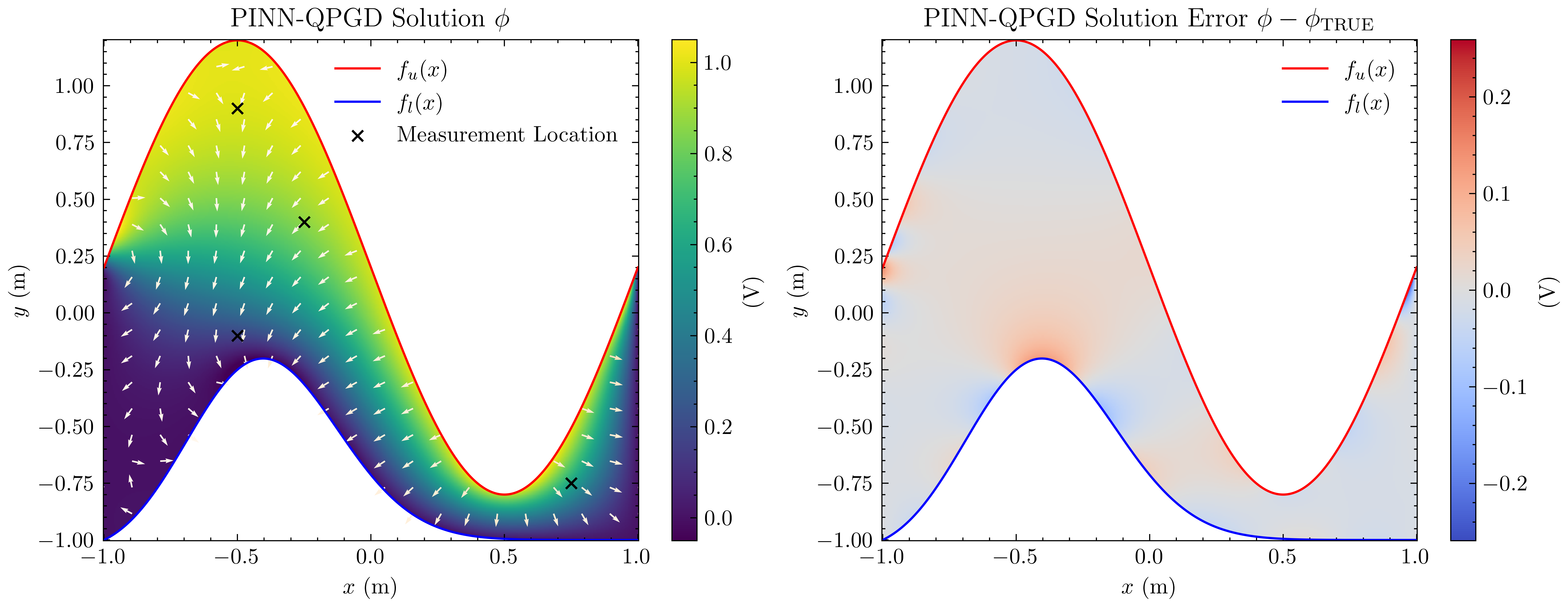}
  \caption{
  \textit{Left}: QPGD solution with $c=1$. Arrows illustrate the normalized direction of the predicted electric field. \textit{Right}: QPGD error with $c=1$.}
  \label{fig:PINN_capacitor_solution}
\end{figure}

In this case study we consider a 2D capacitor with perfect conducting plates of a complex shape. In the interior, between the upper and lower plates, the scalar potential $\phi : \Omega \mapsto \mathbb{R}$ satisfies the 2D Laplace's equation:
\begin{equation}
\nabla^2 \phi (x, y) = 0,
\end{equation}
where
\[
\Omega = \{ (x, y) \in \mathbb{R}^2 : x \in [-1, 1], \, f_{l}(x) \leq y \leq f_{u}(x) \}.
\]
The sides and bottom boundaries are grounded with potential $\phi=0~V$, while the top boundary is at some voltage $\phi=V_0$, see Fig. \ref{fig:PINN_capacitor_solution}:

\begin{align}
\phi(x,f_{u}(x)) &=V_0, \quad x \in [-1, 1], \label{eqn:bc_upper}\\
\phi(0,y) &= 0~V,  \quad y \in (-1, 0.2), \label{eqn:bc_left}\\
\phi(1,y) &= 0~V,  \quad y \in (-1, 0.2), \label{eqn:bc_right}\\
\phi(x,f_{l}(x) )&= 0~V,  \quad x \in [-1, 1], \label{eqn:bc_lower}
\end{align}
where $f_{u}(x) = -\sin (\pi x) + 0.2$ and $f_{l}(x) = e^{-(x+0.5)^2/0.2} (1-x^2) -1 $. 

If $V_0 = 1~V$ is known, the above can be solved with the traditional PINN approach using $l(\param) = l_{\text{PDE}}(\param)+ l_{\text{BC}}(\param)$, for some parameter weights $\param$ -- we use this method to obtain a ground truth solution $\phi_{\text{TRUE}}(x, y) $. 

We consider $V_0$ unknown, with some measurement data of the interior subject to some uncertainty.  The goal is to estimate the $\phi_{\text{TRUE}}(x,y)$ field using a network $\phi(x,y;\param)$ and determine the unknown voltage $V_0$, known as the ``inverse'' problem. The measurement data is corrupted by the addition of Gaussian noise: 
\begin{equation}
\phi_{\text{LABEL}}(x,y) = \phi_{\text{TRUE}}(x, y) +\epsilon, \quad \epsilon \sim \mathcal{N}(0, \delta^2),
\end{equation}
where $\epsilon$ is identically distributed and independent for each measurement, while the uncertainty scale is known: $\delta = 0.1~V$. In real world problems, the uncertainty scale can often be reasonably estimated. 

We construct the losses with the data loss with a $p$-norm:
\begin{align}
l_{\text{DATA}}(\param) &= \left( \frac{1}{N_{M}-1} \sum_{i=1}^{N_{M}} |\phi(x_i, y_i; \param) - \phi_{\text{LABEL}, i}(x_i,y_i)|^p \right)^{1/p} \text{ for } (x_i, y_i) \in \mathcal{M}, \label{eqn:p_norm_loss}
\end{align}
as well as:
\begin{align}
l_{\text{PDE}}(\theta) &= \frac{1}{N_{\text{PDE}}} \sum_{i=1}^{N_{\text{PDE}}} |\nabla^2 \phi(x_i,y_i;\param) |^2 \text{ for } (x_i, y_i) \in \mathrm{int}(\Omega), \label{eqn:pde_loss}\\
l_{\text{BC}, 0}(\param) &=  \frac{1}{N_{\text{BC}, 0}} \sum_{i=1}^{N_{\text{BC},0}} |\phi(x_i, y_i; \param) - 0 |^2 \text{ for } (x_i, y_i) \in \partial \Omega, \label{eqn:bc_lower_loss} \\
l_{\text{BC, V}}(\param) &=  \frac{1}{N_{\text{BC, V}}} \sum_{i=1}^{N_{\text{BC, V}}} |\phi(x_i, y_i; \param) - \hat{V}_0 |^2 \text{ for } (x_i, y_i) \in \partial \Omega, \label{eqn:bc_top_loss}
\end{align} 
where $N_{M}=4$ measurement points from the measurement set $\mathcal{M} \subset \Omega$, $N_{\text{PDE}}=20,000$ PDE points, and $N_{\text{BC}, 0}=300$ boundary condition points for the grounded boundaries (corresponding to \eqref{eqn:bc_left}, \eqref{eqn:bc_right}, \eqref{eqn:bc_lower}), and $N_{\text{BC, V}} = 100$ for upper boundary (corresponding to \eqref{eqn:bc_upper}). The scalar $\hat{V}_0$ is a trainable/learnable parameter, considered as an element of the vector $\theta$ which does not contribute to the output of the network $\phi(x,y; \param)$. For $p=2$, \eqref{eqn:p_norm_loss} is the root mean squared error -- appropriate for Gaussian errors. For $p =\infty $,  \eqref{eqn:p_norm_loss} is maximum absolute error, which is most appropriate for uniformly distributed errors. We formulate the constrained optimization problem: 
\begin{align}
l(\param) &= l_{\text{PDE}}(\param) +  l_{\text{BC}, 0} (\param)
+  l_{\text{BC,V}}(\param) ,   \label{QPGD-Loss} \\
& \text{subject to:}  \nonumber \\
&l_{\text{DATA}}(\param) \leq z \delta, \label{measurement-constraint}
\end{align}
where $z$ is some hyperparameter that quantifies the tolerance for error. For the case of heteroscedastic uncertainties with variable $\delta_i$, (\ref{measurement-constraint}) can be easily modified: $\left( \frac{1}{N_{M}} \sum_i^{N_{M}} |\frac{\phi_i - \phi_{\text{LABEL},i}}{\delta _i}|^p \right)^{1/p} \leq z$. To implement \eqref{eqn:update_law} we take $f(\param) = l(\param)$ and $g(\theta) = l_{\text{DATA}}^2(\param) - z^2 \delta^2$, using squared version of  (\ref{measurement-constraint}) since $l_{\text{DATA}}(\param)^2$ is a scaled mean squared error, similar in form to the mean squared errors in \eqref{eqn:pde_loss} - \eqref{eqn:bc_top_loss}.

\begin{figure}[t]
\centering  \includegraphics[width=0.98\columnwidth]{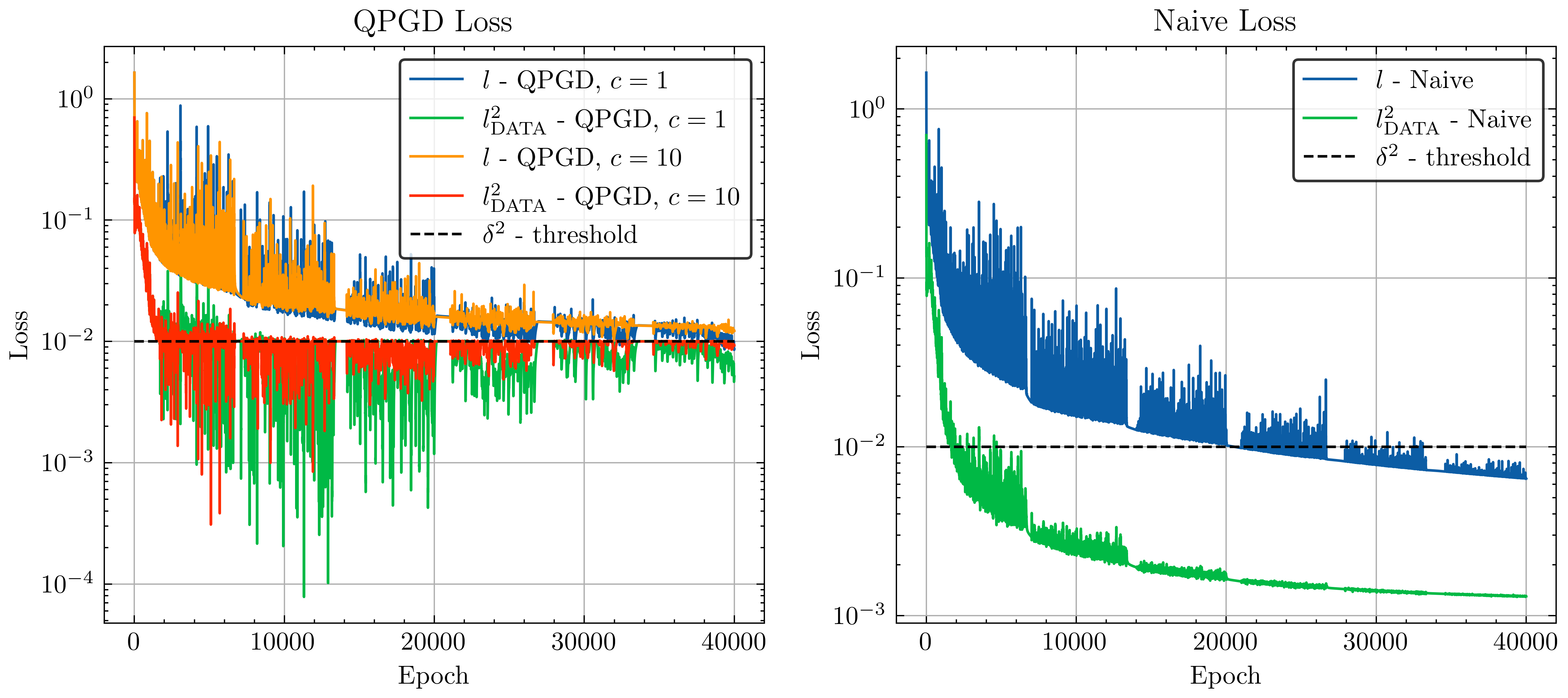}
  \caption{Losses over training for QPGD-PINN and Naive-PINN.}
  \label{fig:loss_plots}
\end{figure}

\begin{figure}[t]
\centering  \includegraphics[width=0.5\columnwidth]{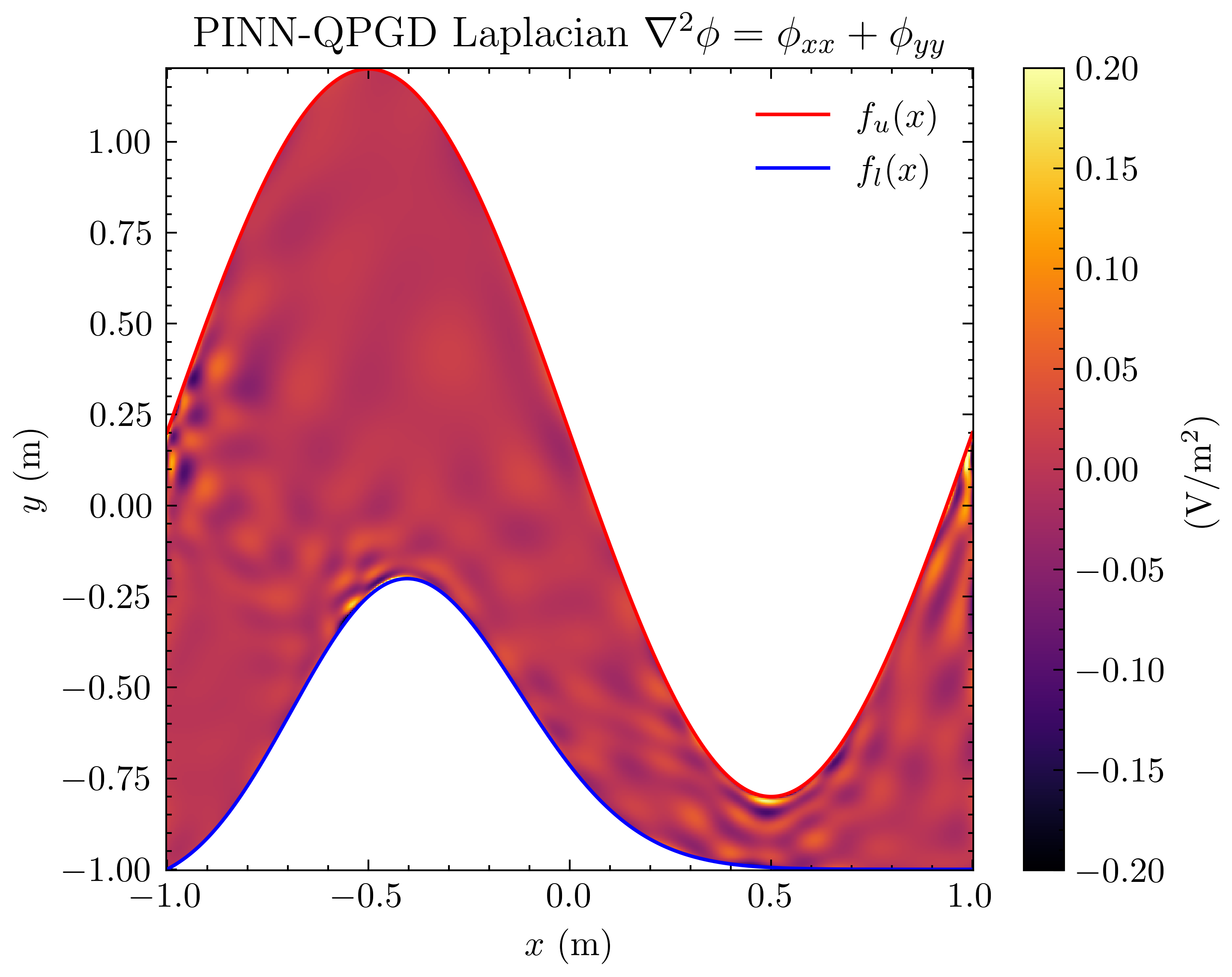}
  \caption{Laplacian $\nabla^2 \phi (x,y;\param)$ plotted over $\Omega$ of the QPGD solution with $c=1$.}
  \label{fig:PINN_Laplacian}
\end{figure}

To test our constraint optimization, we trained a regular PINN (Naive-PINN) with loss $l_{\text{PDE}}(\param) +  l_{\text{BC}, 0} (\param)
+  l_{\text{BC,V}}(\param) + l_{\text{DATA}}(\param)^2$ against two networks 
trained under (\ref{QPGD-Loss}) and (\ref{measurement-constraint}) for $c=1$ and $c=10$, with the latter networks labeled as PINN-QPGD. Since we are given the variance of the noise,  we should expect the fit to follow:
\begin{equation}
\frac{1}{N_{M}-1} \sum_{i=1}^{N_{M}} |\phi(x_i, y_i; \param) - \phi_{\text{LABEL}, i}(x_i,y_i)|^2 \lessapprox \delta^2.    \label{eq:est_variance}
\end{equation} 
We therefore choose $p=2$ and $z=1$. When the left hand side of \ref{eq:est_variance} is less than $\delta^2$, then the network focuses on minimizing $l (\param)$ and avoids overfitting to the data. For each network, we used a dense neural network with 4 hidden layers and 64 neurons per layer, along with the learnable parameter $\hat{V}_0$. For the activation function, we used GELU (Gaussian Error Linear Unit) (\cite{hendrycks2016gaussian}):
$
    \text{GELU}(x) = x \Phi (x)
$
where $\Phi(x) = \frac{1}{\sqrt{2\pi}} \int_{-\infty}^x e^{-t^2/2} dt$. We found this trained better than the usual tanh function used in PINNs, and performed better than a sigmoid as well. This is likely due to the fact GELU, unlike the other two activation functions, does not suffer from the vanishing gradient problem. The GELU curve resembles the commonly used RELU in machine learning, but has the advantage that it is smooth. 

\begin{table}[]
\centering
\begin{tabular}{l|l|l|l}
                           & PINN-Naive  & PINN-QPGD $c=1$   & PINN-QPGD $c=10$ \\ \hline
$\hat{V}_{0}$              & $1.035$ V  & $\textbf{0.996 }$ V  & $0.966$ V   \\ \hline
Average Absolute Error (Interior) & $0.029$ V   & $\textbf{0.017 }$ V & $0.019$ V    \\ \hline
Average Absolute Laplacian & $0.0120$ V/m$^2$ & $0.0114$ V/m$^2$ & $\textbf{0.0106}$ V/m$^2$ \\ \hline
\end{tabular}
\caption{Results of training. Bold indicates best performance.}
\label{Capacitor-Results}
\end{table}

All models were trained on identical learning rate schedules, and initialized with the same weights. The optimizer was one commonly used in machine learning, Adam (\cite{kingma2014adam}), which uses gradient descent along with momentum and RMSProp. Both models were trained for 40,000 epochs, with the learning rate starting at $\gamma = 4\times 10^{-3}$ and decreased by half approximately every 6,6667 epochs. The QPGD-PINN models were pretrained with the Naive-PINN loss until $l_{\text{DATA}^2}$ reached the value $\delta^2$, around epoch 100.

It is possible to train a QPGD-PINN using the update law \eqref{eqn:update_law} without pretraining. We found during our tests that training overall took longer without pretraining. This is likely due to the fact that the set $\{\param: g(\param) \leq 0\}$ likely is not connected, and there may be many ``islands'' that the parameter $\param$ may be driven to due during training -- this certainly violates at least Assumptions \ref{assum:unique_solution} and \ref{assum:gradient_condition} in practice. Therefore we recommend pretraining with a naively designed PINNs loss until the threshold value of the constraint is reached.

The potential field predicted by PINN-QPGD and its error can be seen in Fig. \ref{fig:PINN_capacitor_solution}, while how well it satisfies Laplace's equation can be seen in Fig. \ref{fig:PINN_Laplacian}. As seen in Table \ref{Capacitor-Results}, the PINN-QPGD outperformed PINN-Naive according to several metrics. The PINN-Naive model overfit the data at the price of sacrificing the known physics of the PDE and boundary conditions. With the QPGD scheme, the network can refine the the physics based loss indefinitely after satisfying the data based loss.

The same methodology was also tested using 4 different measurement points, with similar results. When the number of measurement points was increased to 43, naive PINN and PINN-QPGD were about equally effective. This suggest the PINN-QPGD is most useful when the number of data points is low, and there is a higher potential for overfitting. 

\section{Conclusion}
This study presents a methodology for training PINNs for solving forward and inverse problems, addressing critical challenges that have hindered their application. By developing a training scheme based on constrained optimization and a QP-based gradient descent law, QPGD, we provide a solution that simplifies the loss function design process and ensures optimal parameter convergence. This approach balances the competing demands of data-based and PDE residual losses, thereby enhancing the performance and applicability of PINNs in solving physics-based problems. We present a case study of solving Laplace's equation and determining the unknown voltage on a capacitor surface, using noisy measurements of the field in between the plates -- a situation in which overfitting to noise may lead to a poor estimate of the unknown voltage and estimation of the solution of the potential field. Further work can close the gap between theory and  practice, in particular how the inclusion of momentum and RMSProp to gradient descent or the existence of stationary points affects the analysis. Future work will also expand this method by including multiple constraints in the formulation.

\acks{This work was supported by the Los Alamos National Laboratory (LANL) LDRD Program Directed Research (DR) project 20220074DR.}

\section*{Appendix}
Proofs omitted due to space.
\begin{proposition}[Lyapunov Function] \label{prop:unique_min}
Let Assumptions \ref{assum:unique_solution}, \ref{assum:fg_assums}, and \ref{assum:safe_set} hold. Consider a set $\mathcal{S}_{\bar g} = \{g(\param) \leq \bar g \}$ for some $\bar g > 0$. There exists a $\beta^*>0$ such that for all $\beta > \beta^*$ the function $\Vb$ in \eqref{eqn:Vb_definition} has a unique minimizer $\param^*$ on $\mathcal{S}_{\bar g}$.
\end{proposition}

\begin{proposition}[Existence of Small $\epsilon_\alpha$] \label{prop:small_eps}
Let Assumptions \ref{assum:fg_assums}, \ref{assum:safe_set}, and \ref{assum:gradient_condition} hold. There exists an $\epsilon_\alpha^*$ such that for all $\epsilon_\alpha \in (0, \epsilon_\alpha^*)$, $\| \nabla g(\param) \| ^2 \leq \epsilon_\alpha \implies \alpha(\param) = 0$ in \eqref{eqn:alpha}.
\end{proposition}

\begin{proposition}[Forward Invariance of Positive Levels of $g$] \label{prop:forward_invariance}
    Let Assumptions \ref{assum:fg_assums}, \ref{assum:safe_set}, and \ref{assum:gradient_condition} hold and let $\epsilon_\alpha \in (0, \epsilon_\alpha^*)$ for $\epsilon_\alpha^*$ given in Proposition \ref{prop:small_eps}. For any $\bar g > 0$ there exists a $\lr^*$ such that the set $\{g(\param) \leq \bar g \}$ is forward invariant and compact for all $\lr \in (0, \lr^*)$, under the update law \eqref{eqn:update_law}. 
\end{proposition}

\textbf{Proof Sketch of Main Result:}
We first state two preliminaries: 1) $\Vb$ is locally Lipschitz everywhere (it is composed of Lipschitz functions) with a unique minimum of $0$ at $\param^*$ because of Lemma \ref{prop:unique_min}, and $\Vb$ is radially unbounded (Assumption \ref{assum:safe_set}). 2) The right hand side of the update law \eqref{eqn:update_law} is locally Lipschitz with respect to $\paramt$ given Assumptions \assums.

Next we choose $\epsilon_\alpha, \bar g, \beta$. Let $\epsilon_\alpha$ be sufficiently small such that $\| g(\paramt) \|^2 \leq \epsilon_\alpha \implies \alpha( \paramt ) = 0$, which exists by Proposition \ref{prop:small_eps}. This allows cases of $\| g(\paramt) \|^2 \leq \epsilon_\alpha$, in the arguments below, to be a subset of the case of $\alpha(\paramt) = 0$, simplifying the proof and presentation. Choose $\bar g>0$ such that the initial condition $\param^{(0)} \in \{g(\param) \leq \bar g \}$, an invariant set by Proposition \ref{prop:forward_invariance}. Then choose $\beta$ sufficiently large such that i) $\Vb$ has a unique minimum of $0$ at $\param^*$, by Proposition \eqref{prop:unique_min}, ii) $\beta \geq \frac{F^*}{l_h} + \frac{c \bar g}{2 l_h^2} + \frac{1}{2}$ where $F^*$ is the supremum of $\| \nabla f(\param) \|$ on the compact set  $\{ \bar g \geq g(\param) \geq 0 \}$.

 We present the inequalities below, following from expansions and Lipschitzness of the gradients. The points $\param^+$ and $\param$ are shorthand for $\param^{(t+1)}$ and $\param^{(t)}$ respectively. Substituting the update law $\param^+ - \param =- \lr(\nabla f(\param) + \alpha(\param) \nabla g(\param))$:
\begin{align}
    f(\param^+) - f(\param) &\leq \nabla f(\param)^T (\param^+ - \param) + \frac{1}{2} L_f || \param^+ - \param ||^2 \nonumber \\
    & \leq -\lr \nabla f(\param)^T (\nabla f(\param) + \alpha(\param) \nabla g(\param)) + \frac{1}{2} L_f \lr^2 || \nabla f(\param) + \alpha(\param) \nabla g(\param) ||^2 \label{eqn:f_expansion}\\
    g(\param^+) - g(\param) &\leq \nabla g(\param)^T (\param^+ - \param) + \frac{1}{2} L_g || \param^+ - \param ||^2 \nonumber \\
    &\leq -\lr c g(\param) + \frac{1}{2} L_g \lr^2 || \nabla f(\param) + \alpha(\param) \nabla g(\param) ||^2  \label{eqn:g_expansion}
\end{align}
where we use the fact $z - \max\{z, 0\} \leq 0$ and $\| g(\param) \|^2 \leq \epsilon_\alpha \implies \alpha( \param ) = 0$.

Next we show that $\Vb$ is strictly decreasing on the invariant set $\{g(\param) \leq \bar g \}$. Consider the change in $\Vb$ under the update law:
\begin{equation}
    \Delta \Vb = \Vb(\param^+) - \Vb(\param) = f(\param^+) - f(\param) + \beta \max\{g(\param^+), 0 \} - \beta \max\{g(\param), 0 \}.
\end{equation}
We take two main cases, $\alpha(\param) = 0$ for \textbf{case 1}, and $\alpha(\param) > 0$ for \textbf{case 2}, with each of the two cases consisting of two sub-cases involving the signs of $g(\param^+)$. Each case makes use of the expansions \eqref{eqn:f_expansion} - \eqref{eqn:g_expansion}. We summarize the results in each case.

\textbf{Case 1i:} $g(\param^+) \leq 0$. Then with $\lr \leq 1/L_f$:
\begin{equation}
    \Delta \Vb \leq -\frac{1}{2} \lr \| \nabla f(\param) \|^2- \beta \max \{g(\param), 0 \}.
\end{equation}

\textbf{Case 1ii:} $g(\param^+) > 0$ then $\Delta \Vb = f(\param^+) - f(\param) + \beta g(\param^+)- \beta \max \{g(\param), 0 \}$ and
\begin{equation}
    \Delta \Vb \leq -\frac{1}{2} \lr \| \nabla f(\param) \|^2 - \beta c \lr | g(\param) | 
\end{equation}
with $\lr \leq 1/(L_f + \beta L_g)$ and $\lr c \leq 1/2$.

\textbf{Case 2i:} $g(\param^+) \leq 0$. If $g(\param)<0$, then with the choice $\lr \leq 1/L_f$,
\begin{equation}
    \Delta \Vb \leq -\frac{1}{2} \lr u(\param) - \frac{1}{2} \lr \frac{(c g(\param))^2}{\|\nabla g(\param) \|^2}.
\end{equation}
with $u(\param) = \nabla f(\param)^T M(\param) \nabla f(\param)\geq 0$ due to the matrix $M(\param) = I -  \frac{\nabla g(\param) \nabla g(\param)^T}{|| \nabla g(\param)||^{2} }  \succeq 0$. 

If $g(\param) \geq 0$, with $\lr c \leq 1/2$ and $\lr L_f \leq 1$,
\begin{equation}
    \Delta \Vb \leq -\lr \frac{1}{2} u(\param) - \frac{1}{2} \lr c |g(\param)|.
\end{equation}

\textbf{Case 2ii:} $g(\param^+) > 0$.
For $g(\param) < 0$ with $t\leq 1/L_b$, $L_b = L_f + \beta L_g$, 
\begin{equation}
    \Delta \Vb \leq -\lr \frac{1}{2} u(\param) - \frac{1}{2} \lr \frac{(c g(\param))^2}{\|\nabla g(\param) \|^2} - \beta \lr c | g(\param) |.
\end{equation}

For $\bar g \geq g(\param) \geq 0$, and with $\lr \leq 1/L_b$,
\begin{equation}
    \Delta \Vb \leq -\lr \frac{1}{2} u(\param) - \frac{1}{2} \lr c |g(\param)|.
\end{equation}

Assumption \ref{assum:unique_solution} and KKT theory (Proposition 3.3.1, 3.3.4 (\cite{bertsekas1997nonlinear})), guarantees $\Delta \Vb < 0$ everywhere except $\param = \param^*$. Therefore the update law asymptotically converges to $\param^*$ (\cite{bof2018lyapunov}).

\bibliography{main}

\begin{thebibliography}{25}
\providecommand{\natexlab}[1]{#1}
\providecommand{\url}[1]{\texttt{#1}}
\expandafter\ifx\csname urlstyle\endcsname\relax
  \providecommand{\doi}[1]{doi: #1}\else
  \providecommand{\doi}{doi: \begingroup \urlstyle{rm}\Url}\fi

\bibitem[Allibhoy and Cort{\'e}s(2023)]{allibhoy2023control}
Ahmed Allibhoy and Jorge Cort{\'e}s.
\newblock Control barrier function-based design of gradient flows for
  constrained nonlinear programming.
\newblock \emph{IEEE Transactions on Automatic Control}, 2023.

\bibitem[Ames et~al.(2016)Ames, Xu, Grizzle, and Tabuada]{ames2016control}
Aaron~D Ames, Xiangru Xu, Jessy~W Grizzle, and Paulo Tabuada.
\newblock Control barrier function based quadratic programs for safety critical
  systems.
\newblock \emph{IEEE Transactions on Automatic Control}, 62\penalty0
  (8):\penalty0 3861--3876, 2016.

\bibitem[Antonion et~al.(2024)Antonion, Wang, Raissi, and
  Joshie]{antonion2024machine}
Klapa Antonion, Xiao Wang, Maziar Raissi, and Laurn Joshie.
\newblock Machine learning through physics--informed neural networks: Progress
  and challenges.
\newblock \emph{Academic Journal of Science and Technology}, 9\penalty0
  (1):\penalty0 46--49, 2024.

\bibitem[Bertsekas(1997)]{bertsekas1997nonlinear}
Dimitri Bertsekas.
\newblock Nonlinear programming.
\newblock \emph{Journal of the Operational Research Society}, 48\penalty0
  (3):\penalty0 334--334, 1997.

\bibitem[Bischof and Kraus(2021)]{bischof2021multi}
Rafael Bischof and Michael Kraus.
\newblock Multi-objective loss balancing for physics-informed deep learning.
\newblock \emph{arXiv preprint arXiv:2110.09813}, 2021.

\bibitem[Bof et~al.(2018)Bof, Carli, and Schenato]{bof2018lyapunov}
Nicoletta Bof, Ruggero Carli, and Luca Schenato.
\newblock Lyapunov theory for discrete time systems.
\newblock \emph{arXiv preprint arXiv:1809.05289}, 2018.

\bibitem[Cai et~al.(2021)Cai, Mao, Wang, Yin, and Karniadakis]{cai2021physics}
Shengze Cai, Zhiping Mao, Zhicheng Wang, Minglang Yin, and George~Em
  Karniadakis.
\newblock Physics-informed neural networks (pinns) for fluid mechanics: A
  review.
\newblock \emph{Acta Mechanica Sinica}, 37\penalty0 (12):\penalty0 1727--1738,
  2021.

\bibitem[Chen et~al.(2018)Chen, Badrinarayanan, Lee, and
  Rabinovich]{chen2018gradnorm}
Zhao Chen, Vijay Badrinarayanan, Chen-Yu Lee, and Andrew Rabinovich.
\newblock Gradnorm: Gradient normalization for adaptive loss balancing in deep
  multitask networks.
\newblock In \emph{International conference on machine learning}, pages
  794--803. PMLR, 2018.

\bibitem[Cuomo et~al.(2022)Cuomo, Di~Cola, Giampaolo, Rozza, Raissi, and
  Piccialli]{cuomo2022scientific}
Salvatore Cuomo, Vincenzo~Schiano Di~Cola, Fabio Giampaolo, Gianluigi Rozza,
  Maziar Raissi, and Francesco Piccialli.
\newblock Scientific machine learning through physics--informed neural
  networks: Where we are and what’s next.
\newblock \emph{Journal of Scientific Computing}, 92\penalty0 (3):\penalty0 88,
  2022.

\bibitem[Garcia-Cardona and Scheinker(2024)]{garcia2024machine}
Cristina Garcia-Cardona and Alexander Scheinker.
\newblock Machine learning surrogate for charged particle beam dynamics with
  space charge based on a recurrent neural network with aleatoric uncertainty.
\newblock \emph{Physical Review Accelerators and Beams}, 27\penalty0
  (2):\penalty0 024601, 2024.

\bibitem[Hendrycks and Gimpel(2016)]{hendrycks2016gaussian}
Dan Hendrycks and Kevin Gimpel.
\newblock Gaussian error linear units (gelus).
\newblock \emph{arXiv preprint arXiv:1606.08415}, 2016.

\bibitem[Jarolim et~al.(2023)Jarolim, Thalmann, Veronig, and
  Podladchikova]{jarolim2023probing}
Robert Jarolim, JK~Thalmann, AM~Veronig, and Tatiana Podladchikova.
\newblock Probing the solar coronal magnetic field with physics-informed neural
  networks.
\newblock \emph{Nature Astronomy}, 7\penalty0 (10):\penalty0 1171--1179, 2023.

\bibitem[Kingma(2014)]{kingma2014adam}
Diederik~P Kingma.
\newblock Adam: A method for stochastic optimization.
\newblock \emph{arXiv preprint arXiv:1412.6980}, 2014.

\bibitem[Psaros et~al.(2023)Psaros, Meng, Zou, Guo, and
  Karniadakis]{psaros2023uncertainty}
Apostolos~F Psaros, Xuhui Meng, Zongren Zou, Ling Guo, and George~Em
  Karniadakis.
\newblock Uncertainty quantification in scientific machine learning: Methods,
  metrics, and comparisons.
\newblock \emph{Journal of Computational Physics}, 477:\penalty0 111902, 2023.

\bibitem[Raissi et~al.(2019)Raissi, Perdikaris, and
  Karniadakis]{raissi2019physics}
Maziar Raissi, Paris Perdikaris, and George~E Karniadakis.
\newblock Physics-informed neural networks: A deep learning framework for
  solving forward and inverse problems involving nonlinear partial differential
  equations.
\newblock \emph{Journal of Computational physics}, 378:\penalty0 686--707,
  2019.

\bibitem[Raissi et~al.(2020)Raissi, Yazdani, and Karniadakis]{raissi2020hidden}
Maziar Raissi, Alireza Yazdani, and George~Em Karniadakis.
\newblock Hidden fluid mechanics: Learning velocity and pressure fields from
  flow visualizations.
\newblock \emph{Science}, 367\penalty0 (6481):\penalty0 1026--1030, 2020.

\bibitem[Rautela et~al.(2024)Rautela, Williams, and Scheinker]{rautela2024time}
Mahindra Rautela, Alan Williams, and Alexander Scheinker.
\newblock Time-inversion of spatiotemporal beam dynamics using
  uncertainty-aware latent evolution reversal.
\newblock \emph{arXiv preprint arXiv:2408.07847}, 2024.

\bibitem[Scheinker and Pokharel(2023)]{scheinker2023physics}
Alexander Scheinker and Reeju Pokharel.
\newblock Physics-constrained 3d convolutional neural networks for
  electrodynamics.
\newblock \emph{APL Machine Learning}, 1\penalty0 (2), 2023.

\bibitem[Wang et~al.(2021)Wang, Teng, and Perdikaris]{wang2021understanding}
Sifan Wang, Yujun Teng, and Paris Perdikaris.
\newblock Understanding and mitigating gradient flow pathologies in
  physics-informed neural networks.
\newblock \emph{SIAM Journal on Scientific Computing}, 43\penalty0
  (5):\penalty0 A3055--A3081, 2021.

\bibitem[Williams et~al.(2022)Williams, Krsti{\'c}, and
  Scheinker]{williams2022practically}
Alan Williams, Miroslav Krsti{\'c}, and Alexander Scheinker.
\newblock Practically safe extremum seeking.
\newblock In \emph{2022 IEEE 61st Conference on Decision and Control (CDC)},
  pages 1993--1998. IEEE, 2022.

\bibitem[Williams et~al.(2023{\natexlab{a}})Williams, Krsti{\'c}, and
  Scheinker]{williamsCDC2023}
Alan Williams, Miroslav Krsti{\'c}, and Alexander Scheinker.
\newblock Semi-global practical extremum seeking with practical safety.
\newblock In \emph{2023 IEEE 62st Conference on Decision and Control (CDC)}.
  IEEE, 2023{\natexlab{a}}.

\bibitem[Williams et~al.(2023{\natexlab{b}})Williams, Scheinker, Huang, Taylor,
  and Krstic]{williams2023experimental}
Alan Williams, Alexander Scheinker, En-Chuan Huang, Charles Taylor, and
  Miroslav Krstic.
\newblock Experimental safe extremum seeking for accelerators.
\newblock \emph{arXiv preprint arXiv:2308.15584}, 2023{\natexlab{b}}.

\bibitem[Williams et~al.(2024)Williams, Krstic, and
  Scheinker]{williams2024semiglobal}
Alan Williams, Miroslav Krstic, and Alexander Scheinker.
\newblock Semiglobal safety-filtered extremum seeking with unknown cbfs.
\newblock \emph{IEEE Transactions on Automatic Control}, 2024.

\bibitem[Yang and Perdikaris(2019)]{yang2019adversarial}
Yibo Yang and Paris Perdikaris.
\newblock Adversarial uncertainty quantification in physics-informed neural
  networks.
\newblock \emph{Journal of Computational Physics}, 394:\penalty0 136--152,
  2019.

\bibitem[Zhang et~al.(2019)Zhang, Lu, Guo, and
  Karniadakis]{zhang2019quantifying}
Dongkun Zhang, Lu~Lu, Ling Guo, and George~Em Karniadakis.
\newblock Quantifying total uncertainty in physics-informed neural networks for
  solving forward and inverse stochastic problems.
\newblock \emph{Journal of Computational Physics}, 397:\penalty0 108850, 2019.

\end{thebibliography}

\end{document}